\documentclass[11pt,leqno]{article} 
\usepackage{graphics}

\usepackage{lipsum}

\let\OLDthebibliography\thebibliography
\renewcommand\thebibliography[1]{
  \OLDthebibliography{#1}
  \setlength{\parskip}{2pt}
  \setlength{\itemsep}{2pt plus 0.3ex}
}

\newtheorem{thm}{Theorem}[section]
\newtheorem{lma}{Lemma}[section]

\newcommand{\beqa}{\begin{eqnarray}}
\newcommand{\eeqa}{\end{eqnarray}}

\newcommand{\pf}{\noindent {\bf Proof:} $\s$ }
\newcommand{\epf}{ \hfill$\diamondsuit$ \medskip}

\newcommand{\ds}{\displaystyle}
\newcommand{\beq}{\begin{equation}}
\newcommand{\eeq}{\end{equation}}
\newcommand{\lbl}{\label}
\newcommand{\s}{\; \;}

\newcommand{\ep}{\epsilon}

\newcommand{\la}{\lambda}

\newcommand{\ra}{\rightarrow}
\newcommand{\al}{\alpha}
\newcommand{\p}{\varphi}

\title{Infinitely many solutions for a class of resonant problems}

\author{
Philip Korman   \\ 
Department of Mathematical Sciences \\ 
University of Cincinnati \\ 
Cincinnati Ohio 45221-0025 \\
}

\date{}

\begin{document}

\maketitle
\begin{abstract}
\noindent 
We consider radially symmetric solutions for a class of resonant problems on a unit ball $B \subset R^n$ around the origin
\[
\Delta u+\la _1 u +g(u)=f(r)  \s \mbox{for $x \in B$}, \s u=0 \s \mbox{on $\partial B$} \,.
\]
Here the function $g(u)$ is periodic of mean zero, $x \in R^n$, $r=|x|$, $\la _1$ is the principal eigenvalue of $\Delta$ on $B$. The problem has either infinitely many or finitely many solutions depending on the space dimension $n$. The situation turns out to be  different for each of the following cases: $1 \leq n \leq  3$, $n=4$, $n=5$, $n=6$, and $n \geq 7$.
 \end{abstract}

\begin{flushleft}
\end{flushleft}

\begin{flushleft}
\end{flushleft}

\section{Introduction}
\setcounter{equation}{0}
\setcounter{thm}{0}
\setcounter{lma}{0}

Consider a {\em resonant} semilinear problem
\beq
\lbl{i1}
\s\s\s \Delta u+\la _1 u +g(u)=\mu _1 \p _1(x)+e(x)  \s \mbox{for $x \in D$}, \s u=0 \s \mbox{on $\partial D$} \,,
\eeq
on a bounded domain $D \subset R^n$.
Here $x \in R^n$, and $(\la _1, \p _1(x))$ is the principal eigenpair of the Laplacian on $D$, with zero boundary conditions, $\mu _1 \in R$, $e(x) \in \p_1 ^{\perp}$, where $\p_1 ^{\perp}$ denotes the orthogonal complement of $\p _1(x)$ in $L^2(D)$. Solutions of  (\ref{30}) are also decomposed as $u(x)=\xi _1 \p _1(x)+U(x)$, with $U(x) \in \p_1 ^{\perp}$. This problem is {\em at resonance}, which is one of the basic concepts of science and engineering. The study of resonance for elliptic boundary value problem was originated in the classical papers of E.M. Landesman and A.C. Lazer \cite{L}, and A.C. Lazer  and  D.E. Leach \cite{leach}. There is now a huge literature on this topic. In this paper we use a novel approach, based on the theory of global solution curves that are parameterized by the first harmonic of the solution, that we developed in \cite{K1}. Another new ingredient is the use of the stationary phase method to study oscillating integrals. Our approach is uniquely suitable for numerical computations, that we perform and present below. The {\em Mathematica} program was written jointly with D.S. Schmidt, and it is presented  with detailed explanations in \cite{KS1}.
\medskip

In \cite{K1} we gave conditions (see (\ref{30}), (\ref{2a1}) below) under which the solution set of (\ref{i1}) is exhausted by a single continuous solution curve $(u(x),\mu)(\xi _1)$ with the first harmonic of the solution, $\xi _1$, acting as { \em a global parameter}. Namely, for each $\xi _1 \in R$ there is a unique solution pair  $(\la _1, \p _1(x))$. A section of this curve, $\mu _1=\mu _1(\xi _1)$, governs the multiplicity of solutions. In particular, if $\mu _1(\xi _1)$ has infinitely many roots, then the problem 
\beq
\lbl{i2}
\s\s\s \Delta u+\la _1 u +g(u)=e(x)  \s \mbox{for $x \in D$}, \s u=0 \s \mbox{on $\partial D$} 
\eeq
 has infinitely many solutions.
\medskip

For periodic $g(u)$ of mean zero (like $g(u)=\sin u$) a very detailed result was obtained in R. Schaaf  and K. Schmitt \cite{SS1}. They defined $g_1(u),g_2(u),g_3(u)$ to be the unique periodic functions  of mean zero, so that $ g'_1(u)=g(u)$, $ g'_2(u)=g_1(u)$, $ g'_3(u)=g_2(u)$, and showed that the multiplicity depends on the dimension $n$ as follows. For $1 \leq n \leq 3$, the problem (\ref{i2}) has  infinitely many solutions. For $n=4$, a condition on $g_2(u)$ was given for the existence of infinitely many solutions, and a complementary condition was provided under which the number of solutions is finite. For $n \geq 5$ and $g_2(0) \ne 0$,  the number of solutions of  (\ref{i2}) was proved to be  finite.
Only the case $n \geq 5$ and $g_2(0) = 0$ was left open. (That includes $g(u)=\sin u$, with $g_2(u)=-\sin u$, $g_2(0)=0$.)
\medskip

In this paper we investigate the case  $g_2(0) = 0$ for a special class of problems involving radial solutions on a ball around the origin $D=B=\left\{ x \in R^n \; {\rm with} \; ||x||<1   \right\}$. It turned out that for $n = 5$ and $g_2(0) = 0$ the number of solutions is infinite, while for $n \geq 6$ and $g_2(0) = 0$ the number of solutions depends on $g_3(0)$. (We also rederive the results of  \cite{SS1} for the radial case by a different method.) Unlike \cite{SS1}, we use the stationary phase method to study oscillating integrals. The use of global solution curves described above (rather than continuum of solutions as in  \cite{SS1}) allowed us to conclude the radial symmetry of solutions of (\ref{i1}) on $B$, and also that the problem (\ref{i1}) has no solutions for $|\mu _1|$ large (global solution curves  also provide a basis for numerical computations that we present). Once it is established that solutions are radially symmetric on $B$, we proceed similarly to A. Galstyan et al \cite{G}.
\medskip

In the study of oscillating integrals we perform up to three integrations by parts, depending on the dimension $n$, and then use the stationary phase method. We begin with $g(u)=\sin u$, and then generalize.
\medskip 

More general results, without  requiring $g(u)$ to be periodic,  are obtained in the one-dimensional case, by using geometrical arguments instead of stationary phase method.

\section{Radially symmetric oscillatory integrals}
\setcounter{equation}{0}
\setcounter{thm}{0}
\setcounter{lma}{0}

We study oscillating integrals of the form
\beq
\lbl{1}
I(\xi)=\int_0^1 g(\xi v(r)) \p (r) r^{n-1} \, dr \,,
\eeq
depending on a parameter $\xi$. Various choices of the functions $v(r)$ and $\p (r)$ will be considered, beginning with  $v(r)=\p (r)=\p _1(r)$, the principal eigenfunction of the Laplacian on the unit ball $B$, with $\p _1=0$ on $\partial B$. We assume that $g(u)$ is a periodic function of mean zero, which implies  that $g(u)$ changes sign infinitely many times, and the issue is whether $I(\xi)$  changes sign infinitely many times, as $\xi \ra \infty$. It turns out that the answer depends on the dimension $n$.
\medskip

Depending on the dimension $n$, we shall need to perform up to three integrations by parts for $I(\xi)$. Following \cite{SS1}, define $g_1(u),g_2(u),g_3(u)$ to be the unique periodic functions  of mean zero, such that $ g'_1(u)=g(u)$, $ g'_2(u)=g_1(u)$, $ g'_3(u)=g_2(u)$. (In case $g(u)=\sin u$, $g_1(u)=-\cos u$, $g_2(u)=-\sin u$, $g_3(u)=\cos u$.) Denoting $f_1(r)=\frac{\p (r) r^{n-1}}{v'(r)}$, obtain
\beqa
\lbl{2}
& I(\xi)=\frac{1}{\xi} \int_0^1 \frac{\p (r) r^{n-1}}{v'(r)} d \left(g_1(\xi v(r)) \right)=\frac{1}{\xi} \int_0^1 f_1(r) d \left(g_1(\xi v(r)) \right)\\\nonumber
& =\frac{1}{\xi} f_1(r) g_1(\xi v(r)){\Large |}_{_0}^{^1}-\frac{1}{\xi} \int_0^1 f'_1(r) g_1(\xi v(r)) \, dr \,.\\ \nonumber
\eeqa 
Writing $g_1(\xi v)=\frac{1}{\xi v'} \frac{d}{dr} \left(g_2(\xi v(r)) \right)$, and denoting $f_2(r)=\frac{f'_1(r)}{v'(r)}$, integrate by parts again to get
\beqa
\lbl{3}
& I(\xi)=\frac{1}{\xi} f_1(r) g_1(\xi v(r)){\Large |}_{_0}^{^1}-\frac{1}{\xi^2} f_2(r) g_2(\xi v(r)){\Large |}_{_0}^{^1}  \\\nonumber
& + \frac{1}{\xi^2} \int_0^1 f_2'(r) g_2(\xi v(r)) \, dr \,.
\eeqa
Denoting $f_3(r)=\frac{f'_2(r)}{v'(r)}$, and writing $g_2(\xi v)=\frac{1}{\xi v'} \frac{d}{dr} \left(g_3(\xi v(r)) \right)$,
integrate by parts once more to get
\beqa
\lbl{4}
& I(\xi)=\frac{1}{\xi} f_1(r) g_1(\xi v(r)){\Large |}_{_0}^{^1}-\frac{1}{\xi^2} f_2(r) g_2(\xi v(r)){\Large |}_{_0}^{^1}  \\\nonumber
& +\frac{1}{\xi^3} f_3(r) g_3(\xi v(r)){\Large |}_{_0}^{^1}-\frac{1}{\xi^3} \int_0^1 f_3'(r) g_3(\xi v(r)) \, dr \,.
\eeqa

We shall use the following lemma, based on the stationary phase method, see e.g., A. Galstian et al \cite{G} or P. Korman \cite{K}.
\begin{lma}
\lbl{lma:1}
Assume that the functions  $f(x)$ and  $\p (x)>0$ are of class  $C^2[0,1]$, and satisfy
\[
\p '(x)<0 \quad {\rm for \, all}\quad x\in(0,1], \quad {\rm and }\quad \p '(0)=0,\,\, \p ''(0)<0\,.
\]
Then, as $\xi \to \infty$,
\[
\int_0^1 f(x)e^{i\xi \p (x)}dx =  e^{i(\xi \p (0)-
\frac{\pi}{4})}\sqrt {\frac{\pi}{2\xi |\p ''(0)|}}f(0) + O\left(\frac{1}{\xi}\right)\,.
\]
\end{lma}

The proof of this lemma can be found in e.g., \cite{G}, however the proof is sketched next, since  a similar idea is used later on to prove a more general result. The term $e^{i\xi \p (x)}$ involves fast oscillations about zero, which are mutually cancelling, except near $x=0$, where $\p (x) \approx \p (0)+\frac{1}{2} \p ''(0)x^2$, and the oscillations are slow. Then
\[
\int_0^1 f(x)e^{i\xi \p (x)}dx \approx  e^{i\xi \p (0)} \int_0^1 f(x)e^{i\xi \frac{1}{2} \p ''(0)x^2}dx \,.
\]
The  evaluation of the resulting Frenet-type integral  is contained in the following lemma, see e.g., A. Galstian et al \cite{G} or P. Korman \cite{K}.
\begin{lma}\lbl{lma:2}
Assume that $f(x) \in C^2[0,a]$ for some number $ a >0$. Then as $\xi \ra \infty$
\[
\int_0^a f(x)e^{\frac{1}{2} i\al \xi x^2} \, dx=e^{i \frac{\pi}{4} \delta (\al )} \sqrt{ \frac{\pi}{2|\al| \xi }}f(0)+O \left( \frac{1}{\xi} \right) \,,
\]
where $\delta (\al )={\rm sign } \, \al$.
\end{lma}

Recently the Lemma \ref{lma:1} was used to study oscillatory bifurcation curves by T. Shibata \cite{sh} and K. Kato and T. Shibata \cite{ks}.
\medskip

We start by  considering  a model case
\beq
\lbl{5}
J(\xi)=\int_0^1 \sin \left(\xi \p _1(r) \right) \p _1(r) r^{n-1} \, dr \,,
\eeq
where $\p _1(r)$ is the principal eigenfunction of the Laplacian on the unit ball  $B \subset R^n$, $\p _1(r)=c_0 r^{-\frac{n-2}{2}} J_{\frac{n-2}{2}}(\nu _1 r)$, $\p _1(1)=0$. Here  $\nu _1>0$ denotes the first root of the Bessel function $J_{\frac{n-2}{2}}(r)$, and 
$c_0$ is chosen so that $\p _1(0)=1$. The corresponding principal eigenvalue is $\la _1=\nu _1^2$. From the equation 
\beq
\lbl{5.1}
\p _1''+\frac{n-1}{r}\p _1'+\la _1 \p _1=0
\eeq
it follows that
\beq
\lbl{6}
\p _1''(0)=-\frac{\nu _1^2}{n}=-\frac{\la _1}{n} \,.
\eeq

The following result is similar to that in A. Galstian et al \cite{G}. The case $1 \leq n \leq 4$ was covered previously in R. Schaaf  and K. Schmitt \cite{SS1} by a different method. Our approach provides accurate asymptotic formulas, in addition to the oscillation properties.

\begin{thm}\lbl{thm:12}
For $1 \leq n \leq 5$, $J(\xi)$ changes sign infinitely many times on $(0,\infty)$, while for $n \geq 6$ the number of sign changes is at most finite.
\end{thm}

\pf
We shall use the integration by parts formulas (\ref{2}), (\ref{3}), (\ref{4}) with $v(r)=\p (r)=\p _1(r)$. Here $g(u)=\sin u$, $g_1(u)=-\cos u$, $g_2(u)=-\sin u$, $g_3(u)=\cos u$, and also 
\beq
\lbl{6.5}
f_1(r)=\frac{r^{n-1} \p _1(r)}{\p _1'(r)} \,,  
\eeq
\beq
\lbl{6.5a}
  f_2(r)=\frac{f_1'(r)}{\p _1'(r)}=\frac{r^{n-2} \left[(n-1)\p _1 \p _1'-r \p _1\p _1''+r {\p _1'}^2 \right]}{{\p _1'}^3} \,.
\eeq
Expressing the second and the third derivatives of $\p _1(r)$ from (\ref{5.1}), obtain
\beqa
\lbl{6.5b}
& f_3(r)=\frac{f_2'(r)}{\p _1'(r)}=\frac{r^{n-3}}{{\p _1'}^5(r)} {\large [3\la_1^2 r^2 \p _1^3+8(n-1)\la _1 r \p_1^2 \p'_1} \\
& + \left( 8-14n+6n^2+3 \la _1r^2  \right)\p_1 {\p _1'}^2+4(n-1)r {\p _1'}^3   {\large ]}  \,.          \nonumber
\eeqa

We consider the following cases  depending on the dimension $n$.
\medskip

\noindent
{\bf i.} $\ds n=2$ (the case $n=1$ is similar). Here $f_1(r)=\frac{r \p _1 (r)}{\p' _1 (r)}$, $f_1(1)=0$, $f_1(0)=\frac{\p _1(0)}{\p _1''(0)}=-\frac{2}{\nu _1^2}$. It is straightforward to verify that $f_1(r) \in C^{\infty}[0,1)$ for $n \geq 2$, see \cite{G}. By (\ref{2}) and Lemma \ref{lma:1}
\beqa \nonumber
& J(\xi)= -\frac{1}{\xi} f_1(r) \cos \left(\xi \p _1(r) \right) {\Large |}_{_0}^{^1}  +\frac{1}{\xi} \int_0^1 f'_1(r) \cos \left(\xi \p _1(r) \right) \, dr\\ \nonumber
& =-\frac{2}{\nu _1^2 \xi} \cos \xi  +O \left(\frac{1}{\xi ^{\frac32}} \right) \,,
\eeqa
so that $J(\xi)$ changes sign infinitely many times, as $\xi \ra \infty$.
\medskip

\noindent
{\bf ii.} $\ds n=3$. Now $f_1(r)=\frac{r^2 \p _1 (r)}{\p' _1 (r)}$, $f_1(0)=f_1(1)=0$, while by (\ref{6.5a}) $ f'_1(0)=\frac{ \p _1(0)}{\p ''_1(0)} \ne 0$. Then using (\ref{2}) and Lemma \ref{lma:1} again
\beqa 
\lbl{7}
& J(\xi)= \frac{1}{\xi} \int_0^1 f'_1(r) \cos \left(\xi \p _1(r) \right) \, dr= \frac{1}{\xi} \, {\rm Re}  \int_0^1 f'_1(r) e^{i\xi \p _1(r) } \, dr\\ \nonumber
& =\frac{f_1'(0)}{\xi ^{\frac32}}  \sqrt{\frac{\pi}{2 |\p _1''(0)|}} \cos \left(\xi  -\frac{\pi}{4}\right)+O \left(\frac{1}{\xi ^2} \right) \,,
\eeqa
so that $J(\xi)$ changes sign infinitely many times.
\medskip

\noindent
{\bf iii.} $\ds n=4$. It is straightforward to verify that $f_2(r) \in C^{\infty}[0,1)$ for $n \geq 4$, see \cite{G}. Here $f_1(r)=\frac{r^3 \p _1 (r)}{\p' _1 (r)}$. Again we have  $f_1(0)=f_1(1)=0$, while now $f'_1(0)=0$, so that the principal term in (\ref{7}) is zero. One needs to integrate by parts again, i.e. to use (\ref{3}):
\[
J(\xi)= \frac{1}{\xi ^2} f_2(r) \sin \left(\xi \p _1(r) \right) {\Large |}_{_0}^{^1}  -\frac{1}{\xi ^2} \int_0^1 f'_2(r) \sin \left(\xi \p _1(r) \right) \, dr \,.
\]
The first term is equal to $-\frac{1}{\xi ^2} f_2(0) \sin \xi $, with $f_2(0)=\frac{2}{{\p _1''(0)}^2} \ne 0$  by (\ref{6.5a}), while the integral term is $ O \left(\frac{1}{\xi ^{\frac52}} \right) $ by  Lemma \ref{lma:1}, so that  $J(\xi)$ changes sign infinitely many times.
\medskip

\noindent
{\bf iv.} $\ds n=5$. Now $f_1(r)=\frac{r^4 \p _1 (r)}{\p' _1 (r)}$, $f_1(0)=f_1(1)=0$, and $f_2(0)=0$, while $f'_2(0) \ne 0$  by (\ref{6.5b}) and (\ref{6}) . By (\ref{3}) and Lemma \ref{lma:1}
\beqa \nonumber
& J(\xi)=   -\frac{1}{\xi ^2} \int_0^1 f'_2(r) \sin \left(\xi \p _1(r) \right) \, dr=- \frac{1}{\xi^2 } \, {\rm Im}  \int_0^1 f'_2(r) e^{i\xi \p _1(r) } \, dr\\ \nonumber
& =-\frac{f_2'(0)}{\xi  ^{\frac52}}  \sqrt{\frac{\pi}{2 |\p _1''(0)|}} \sin \left(\xi  -\frac{\pi}{4}\right)+O \left(\frac{1}{\xi ^3} \right) \,,
\eeqa
so that $J(\xi)$ changes sign infinitely many times.
\medskip

\noindent
{\bf v.} $\ds n=6$. Now $f_1(r)=\frac{r^5 \p _1 (r)}{\p' _1 (r)}$, $f_1(0)=f_1(1)=f_2(0)=0$, and also  $f'_2(0)= 0$. We need to integrate by parts one more time, i.e., to use (\ref{4}). It is straightforward to verify that $f_3(r) \in C^{\infty}[0,1)$ for $n \geq 6$, see \cite{G}.
Obtain
\beqa 
\lbl{8}
& J(\xi)= \frac{1}{\xi ^3} f_3(r) \cos \left(\xi \p _1(r) \right) {\Large |}_{_0}^{^1}  -\frac{1}{\xi ^3} \int_0^1 f'_3(r) \cos \left(\xi \p _1(r) \right) \, dr \\ \nonumber\nonumber
& = \frac{1}{\xi ^3} f_3(1)- \frac{1}{\xi ^3} f_3(0) \cos \xi  -\frac{1}{\xi ^3} \int_0^1 f'_3(r) \cos \left(\xi \p _1(r) \right) \, dr \,.\nonumber
\eeqa
The integral term is $O \left(\frac{1}{\xi ^{\frac72}} \right)$ by  Lemma \ref{lma:1}. Whether $J(\xi)$ changes sign finitely or infinitely many times will depend on the relative sizes of $|f_3(1)|$ and $|f_3(0)|$. Using {\em Mathematica}, one calculates $f_3(1) \approx 71.44$ and $f_3(0) \approx -0.09$. It follows that the term $ \frac{1}{\xi ^3} f_3(1)$ is dominant in $J(\xi)$ for large $\xi$, implying that $J(\xi)$ changes sign at most finitely 
 many times.
\medskip

\noindent
{\bf vi.} $\ds n \geq 7$. Now $f_3(0)=0$,  $f_3(1)=\frac{4(n-1)}{{\p _1'}^2(1)} \ne 0$, and by the formula (\ref{8}),  $J(\xi) \sim \frac{1}{\xi ^3} f_3(1)$ for large $\xi$, and hence $J(\xi)$ changes sign at most finitely 
 many times.
\epf 

We turn to the oscillation of more general integrals
\beq
\lbl{10}
K(\xi)=\int_0^1 g \left(\xi \p _1(r) \right) \p _1(r) r^{n-1} \, dr \,,
\eeq
with periodic $g(u)$ of mean zero.  It turns out that for  $n \geq 4$, the number of oscillations will depend on $g(u)$, particularly on whether $g_2(0)$ is zero or not. (The case $g_2(0) \ne 0$ for $n=4$ was already considered in \cite{SS1}).)
We shall need the following generalization of Lemma \ref{lma:1}.
\begin{lma}
\lbl{lma:3}
Assume that   $f(x),\p (x) \in C^2[0,1]$, $h(u) \in C^1(R)$, and $x=0$ is the unique critical point of  $\p (x)$ on $[0,1)$, with  $\p '(x)<0$ on $(0,1)$, and  $\p '(0)=0$, $\p ''(0) < 0$. Assume also that $|h'\left(\xi \p (0) \right)| >0$.
Then, as $\xi \to \infty$,
\beq
\lbl{12}
\s\s\s \int_0^1 f(x)e^{ih\left(\xi \p (x) \right)} \, dx =  e^{i \left[h(\xi \p (0))-\delta
\frac{\pi}{4} \right]}\sqrt {\frac{\pi}{2\xi |h'\left(\xi \p (0) \right)\p ''(0)|}}f(0) + O\left(\frac{1}{\xi}\right)\,,
\eeq
where $\delta= {\rm sign} \,  \left( h'\left(\xi \p (0) \right) \right)$.
\end{lma}

\pf
Let $p(x)=h\left(\xi \p (x) \right)$.  Calculate $p(0)=h(\xi \p (0))$, $p'(x)=\xi h'(\xi \p (x)) \p '(x)$, so that $p'(0)=0$, since $\p '(0)=0$. Then $p''(0)=\xi h'(\xi \p (0)) \p ''(0)$, using again that  $\p '(0)=0$. As above, we approximate $p(x) \approx p(0)+\frac{1}{2}p''(0)x^2=p(0)+\frac{1}{2}\xi h'(\xi \p (0)) \p ''(0)x^2$, for small $x$. Then the integral in (\ref{12}) is approximated by
\[
\int_0^1 f(x)   e^{i \left[h(\xi \p (0))+ \frac{1}{2} \xi h'(\xi \p (0)) \p ''(0) x^2 \right]}   \, dx \,.
\]
Application of the Lemma \ref{lma:2} gives the asymptotic formula above.
As in \cite{G}, the derivation above is justified  by a change of variables $x \ra t$, given by $\p (x)-\p (0)=-t^2$, transforming $\p (x)$ to its quadratic part.
\epf

Taking the imaginary part of (\ref{12}), gives
\beqa
\lbl{14}
& \hspace{.3in} \int_0^1 f(x) \sin h\left(\xi \p (x) \right) \, dx =  \sqrt {\frac{\pi}{2\xi |h'\left(\xi \p (0) \right)\p ''(0)|}}f(0) \sin  \left(h(\xi \p (0)) 
 -\delta
\frac{\pi}{4} \right) \\
& + O\left(\frac{1}{\xi}\right) , \s\s \mbox{with $\delta= {\rm sign} \,  \left( h'\left(\xi \p (0) \right) \right)$} \,. \nonumber
\eeqa

We shall use this formula to study oscillations of the integral $K(\xi )$ in (\ref{10}), by writing $h(u)=\sin ^{-1} g(u)$. For periodic $g(u)$ there are infinitely many points where $h'\left(\xi \p (0) \right)=0$, near which the asymptotic formula (\ref{14}) is not accurate. However, we shall argue that the formula  (\ref{14}) is  accurate on infinitely many intervals where $K(\xi)$ takes both negative and positive values, implying that $K(\xi )$ changes sign  infinitely many times.
\medskip

\noindent
{\bf \large Example } For $g(u)=\sin ^3 u$, $h(u)=\sin ^{-1} \left( g(u) \right)$, and $\p (x)=1-\frac{1}{2}x^2$ we used {\em Mathematica} to calculate both the integral in (\ref{14}) (solid line), and its asymptotic approximation  in (\ref{14}) (dashed line), plotted in Figure \ref{fig:2}. Figure \ref{fig:2} shows that  the asymptotic formula is accurate on sufficiently many intervals to conclude that the integral  changes sign  infinitely many times.

\begin{figure}
\begin{center}
\scalebox{0.99}{\includegraphics{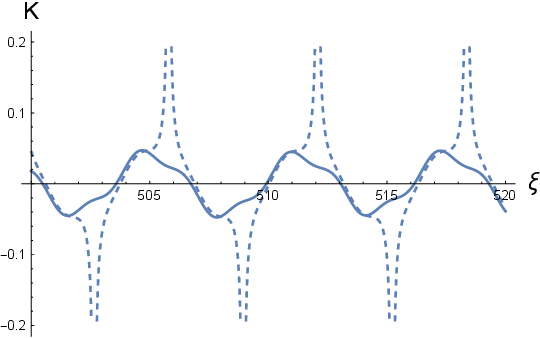}}
\end{center}
\caption{ The  integral in (\ref{14}), and its asymptotic approximation}
\lbl{fig:2}
\end{figure}

\medskip

We shall use the following generalization of  the Riemann-Lebesgue lemma, which is included in Theorem 3 of O. Costin et al \cite{C0}.

\begin{lma}
\lbl{lma:4}
Let $g(u) \in C(0,\infty)$ be a periodic function of mean zero,  and $f(x) \in L^1(0,1)$. Then
\[
\int_0^1 g(\xi x)f(x) \, dx \ra 0 \s \mbox{as $\xi \ra \infty$} \,. 
\]
\end{lma}

In the following theorem we shall apply the asymptotic formula (\ref{14}) to  periodic functions $h(u)$ of the form $h(u)=\sin ^{-1} g(u)$ with various periodic $g(u)$. The formula (\ref{14}) does not apply at the infinitely many roots of $h'(\xi \p (0))$. We shall argue that such points are rare, and the formula  (\ref{14}) does apply at infinitely many points $\xi $, at which $\sin \left( h(\xi \p (0)) \right)$ takes both positive and negative values.
\medskip

\begin{thm}\lbl{thm:2}
Let $g(u) \in C^1(0,\infty)$ be a periodic function of mean zero.
\smallskip

\noindent
(i) For $1 \leq n \leq 3$ the integral $K(\xi)$ (given by (\ref{10})) changes sign infinitely many times on $(0,\infty)$.
 \smallskip

\noindent
(ii) In case $n=4$, $K(\xi)$ changes sign infinitely many times, provided that the same is true for the function $- \frac{1}{\p _1'(1)} g_2(0)+ \frac{2}{{\p _1 ''}^2(0)}  g_2 \left(\xi  \right)$ (which includes the case $g_2(0)=0$), and there are only finitely many sign changes otherwise.
\smallskip

\noindent
(iii) For  $n=5$, in case $g_2(0) \ne 0$ there are only finitely many  sign changes of  $K(\xi)$, and if $g_2(0) = 0$, $K(\xi)$ changes sign infinitely many times. 
\smallskip

\noindent
(iv) For $n = 6$, in case $g_2(0) \ne 0$ there are only finitely many  sign changes of  $K(\xi)$, and if $g_2(0) = 0$,  $K(\xi)$ changes sign infinitely many times, provided that the same is true for the function $ f_3(1)g_3(0)-  f_3(0) g_3\left(\xi  \right)$ (which includes the case $g_3(0)=0$), and there are only finitely many sign changes otherwise.
\smallskip

\noindent
(v) For $n \geq 7$, assume that  $g_2(0) = 0$, but  $g_3(0) \ne 0$. Then the number of sign changes of   $K(\xi)$ is at most finite.
\end{thm}

\pf
The proof is similar to that of Theorem \ref{thm:12}. The functions $f_1(r), f_2(r),f_3(r)$ are the same, as given in (\ref{6.5}),(\ref{6.5a}),(\ref{6.5b}). The breakdown into cases is similar.
\medskip

\noindent
{\bf i.} $\ds n=2$ (the case $n=1$ is similar). As in Theorem \ref{thm:12},  $f_1(r)=\frac{r \p _1 (r)}{\p' _1 (r)}$, $f_1(1)=0$, $f_1(0)=\frac{\p _1(0)}{\p _1''(0)}=-\frac{2}{\nu _1^2}$, and $f_1(r) \in C^{\infty}[0,1)$ for $n \geq 2$. By (\ref{2}) and Lemma \ref{lma:4}
\beqa \nonumber
& K(\xi)= \frac{1}{\xi} f_1(r) g_1 \left(\xi \p _1(r) \right) {\Large |}_{_0}^{^1}  -\frac{1}{\xi} \int_0^1 f'_1(r) g_1 \left(\xi \p _1(r) \right) \, dr\\ \nonumber
& =\frac{2}{\nu _1^2 \xi} g_1 \left(\xi  \right)+o \left(\frac{1}{\xi } \right) \,,
\eeqa
so that $K(\xi)$ changes sign infinitely many times.
\medskip

\noindent
{\bf ii.} $\ds n=3$. As in Theorem \ref{thm:12},  $f_1(r)=\frac{r^2 \p _1 (r)}{\p' _1 (r)}$, $f_1(0)=f_1(1)=0$, while $f'_1(0)=\frac{ \p _1(0)}{\p ''_1(0)} < 0$. By (\ref{2}) 
\[
K(\xi)=-\frac{1}{\xi} \int_0^1 f'_1(r) g_1 \left(\xi \p _1(r) \right) \, dr \,.
\]
Denoting $G_1=\max _{(-\infty,\infty)} |g_1(u)|$, write this integral as $G_1 \int_0^1 f'_1(r) \sin h \left(\xi \p _1(r) \right) \, dr$, with $h(u)=\sin ^{-1} \left( \frac{g_1(u)}{G_1} \right)$ and use (\ref{14}) to express
\beqa
\lbl{16}
& \frac{K(\xi)}{G_1}= \frac{1}{\xi} \sqrt {\frac{\pi}{2\xi |h'\left(\xi  \right)\p _1''(0)|}}f_1'(0) \sin  \left(h(\xi ) 
 -\delta
\frac{\pi}{4} \right)\\ \nonumber
& + O\left(\frac{1}{\xi ^2}\right) \,, \nonumber
\eeqa
with $\delta= {\rm sign} \,  \left( h'(\xi ) \right)= {\rm sign} \,  \left( g_1'(\xi ) \right)$. We now show that $\sin  \left(h(\xi ) 
 -\delta
\frac{\pi}{4} \right)$ changes sign infinitely many times.
The function $\frac{g_1(\xi)}{G_1} $ is periodic in $\xi$, changing  sign infinitely many times. Let $\xi _0$ be a point where $g_1(\xi _0)=0$, and $g_1(\xi)$ changes sign from negative to positive across $\xi _0$. To the right of $\xi _0$ we can find a point $\xi _1$ where $g_1'(\xi _1)>0$, so that $\delta =1$ in (\ref{16}), and $h(\xi _1)$ is small. By (\ref{16}), $K(\xi _1)>0$ (if necessary, adding to $\xi _1 $ a multiple of the period of $g_1(u)$ to make the first term in (\ref{16}) dominant). By the periodicity of $h(\xi  )$ we have a sequence $\{ \xi _n \} \ra \infty$ such that $K(\xi _n)<0$. Similarly, there is  a sequence $\{ \eta _n \} \ra \infty$ such that $K(\eta _n)>0$. Hence,
$K(\xi)$ changes sign infinitely many times.
\medskip

\noindent
{\bf iii.} $\ds n=4$.  Here $f_1(r)=\frac{r^3 \p _1 (r)}{\p' _1 (r)}$, $f_2(r)=\frac{f'_1(r)}{\p' _1 (r)}=\frac{r^2 \varphi _1 '(r) \left(r
   \varphi _1 '(r)+3 \varphi _1
   (r)\right)-r^3 \varphi _1 (r)
   \varphi _1 ''(r)}{\varphi _1 '^{3}(r)}$. Again we have  $f_1(0)=f_1(1)=0$, while now $f'_1(0)=0$, so that the principal term in (\ref{14}) is zero. As in Theorem \ref{thm:12}, $f_2(r) \in C^{\infty}[0,1)$ for $n \geq 4$. One needs to integrate by parts again, i.e., to use (\ref{3}):
\beqa
\lbl{17}
&  K(\xi)= -\frac{1}{\xi ^2} f_2(1) g_2(0)+\frac{1}{\xi ^2}f_2(0)g_2 \left(\xi  \right)  \\ \nonumber
& +\frac{1}{\xi ^2} \int_0^1 f'_2(r) g_2 \left(\xi \p _1(r) \right) \, dr \,. \nonumber
\eeqa
The integral term is $o\left(\frac{1}{\xi ^2} \right)$ by Lemma \ref{lma:4}. Indeed, 
\[
\int_0^1 f'_2(r) g_2 \left(\xi \p _1(r) \right) \, dr=\int_0^1 g_2 \left(\xi x \right) f'_2\left( \p _1^{-1}(x)\right) \psi (x) \, dx \,,
\]
with $\psi (x) \equiv \frac{d}{dx} \p _1^{-1}(x) \in L^1(0,1) $.
Hence, $K(\xi)$ changes sign infinitely many times, provided that the same is true for the function $-\frac{1}{\xi ^2} f_2(1) g_2(0)+\frac{1}{\xi ^2}f_2(0)g_2 \left(\xi  \right)=-\frac{1}{\xi ^2} \frac{1}{\p _1'(1)} g_2(0)+\frac{1}{\xi ^2} \frac{2}{{\p _1 ''}^2(0)}  g_2 \left(\xi  \right)$ (which includes the case $g_2(0)=0$), and there are only finitely many sign changes otherwise.
\medskip

\noindent
{\bf iv.} $\ds n=5$.  Here $f_1(r)=\frac{r^4 \p _1 (r)}{\p' _1 (r)}$, $f_2(r)=\frac{r^3 \varphi _1 '(r) \left(r
   \varphi _1 '(r)+4 \varphi _1
   (r)\right)-r^4 \varphi _1 (r)
   \varphi _1 ''(r)}{\varphi _1 '^{3}(r)}$. In addition to   $f_1(0)=f_1(1)=0$, we now have $f_2(0)=0$. Also  $\ds f'_2(0) =\frac{3 \la _1^2}{25 {\p _1''}^4(0)} \ne  0$, as follows by expressing $f_2'$ from (\ref{6.5b}). The formula (\ref{17}) becomes
\[
K(\xi)= -\frac{1}{\xi ^2} f_2(1) g_2(0)  +\frac{1}{\xi ^2} \int_0^1 f'_2(r) g_2 \left(\xi \p _1(r) \right) \, dr \,.
\]
The integral term is $o\left(\frac{1}{\xi ^2} \right)$ by Lemma \ref{lma:4}. In case $g_2(0) \ne 0$ there are only finitely many  sign changes for $K(\xi )$. Observe that the same is true for all $n \geq 5$. In case  $g_2(0) = 0$, $K(\xi)$ changes sign infinitely many times, using the argument similar to the case $n=3$. 
\medskip

\noindent
{\bf v.} $\ds n=6$. We assume that $g_2(0) = 0$, since in case $g_2(0) \ne 0$ there are only finitely many  sign changes for $K(\xi )$, as we just saw. Now $f_1(r)=\frac{r^5 \p _1 (r)}{\p' _1 (r)}$,  $f_2(r)=\frac{r^3 \varphi _1 '(r) \left(r
   \varphi _1 '(r)+4 \varphi _1
   (r)\right)-r^4 \varphi _1 (r)
   \varphi _1 ''(r)}{\varphi _1 '^{3}(r)}$. Calculate $f_1(0)=f_1(1)=f_2(0)=0$. Also, $f_2'(0)=0$. We need to integrate by parts one more time, i.e., to use (\ref{4}). It is straightforward to verify that $f_3(r) \in C^{\infty}[0,1)$ for $n \geq 6$.
Obtain
\beqa 
\lbl{18}
& K(\xi)= \frac{1}{\xi ^3} f_3(r) g_3\left(\xi \p _1(r) \right) {\Large |}_{_0}^{^1}  -\frac{1}{\xi ^3} \int_0^1 f'_3(r) g_3\left(\xi \p _1(r) \right) \, dr \\ \nonumber\nonumber
& = \frac{1}{\xi ^3} f_3(1)g_3(0)- \frac{1}{\xi ^3} f_3(0) g_3\left(\xi  \right) -\frac{1}{\xi ^3} \int_0^1 f'_3(r) g_3\left(\xi \p _1(r) \right) \, dr \,.\nonumber
\eeqa
As in Theorem \ref{thm:12}, $f_3(1) \approx 71.44$ and $f_3(0) \approx -0.09$. The integral term is $o\left(\frac{1}{\xi ^3} \right)$ by Lemma \ref{lma:4}. $K(\xi)$ changes sign infinitely many times, provided that the same is true for the function $ f_3(1)g_3(0)-  f_3(0) g_3\left(\xi  \right)$ (which includes the case $g_3(0)=0$), and there are only finitely many sign changes otherwise.
\medskip

\noindent
{\bf vi.} $\ds n \geq 7$. Now $f_3(0)=0$,  $f_3(1) \ne 0$, and by the  formula (\ref{18}),  $K(\xi)$ changes sign at most finitely 
 many times.
\epf 

\noindent
{\bf Remark } What if $n \geq 7$, but $g_2(0)=g_3(0)=0$? It appears that infinitely many oscillations are still possible for such special functions $g(u)$, but one would need more than three integrations by parts for a proof.

\section{Oscillations of  the solution curve }
\setcounter{equation}{0}
\setcounter{thm}{0}
\setcounter{lma}{0}

We now consider the following Dirichlet problem on a unit ball $B \subset R^n$ around the origin
\beq
\lbl{30}
\s\s\s \Delta u+\la _1 u +g(u)=f(r)=\mu _1 \p _1(r)+e(r)  \s \mbox{for $x \in B$}, \s u=0 \s \mbox{on $\partial B$} \,.
\eeq
Here $x \in R^n$, $r=|x|$ and $(\la _1, \p _1(r))$ is the principal eigenpair of the Laplacian on $B$, with zero boundary conditions, $\mu _1 \in R$, $e(r) \in \p_1 ^{\perp}$ in $L^2(B)$, and $e(r) \in C^{\al }(B)$, for some $\al \in (0,1)$. Solutions of  (\ref{30}) are decomposed as $u(r)=\xi _1 \p _1(r)+U(r)$, with $U(r) \in \p_1 ^{\perp}$ in $L^2(B)$. The following result describes all solutions of (\ref{30}).
\begin{thm}\lbl{thm:1}
Assume that $g(u) \in C^2(R)$,  and 
\beq
\lbl{2a1} 
g'(u) <\la _2- \la _1 , \s \mbox{for all $u \in R$}\,,
\eeq
\beq
\lbl{2a2}
|g(u)|<\gamma |u|+c , \s \mbox{with $0<\gamma <\la _2- \la _1$, $c \geq 0$, and  $u \in R$} \,.
\eeq
Then the solution set of (\ref{30}) consists of a single continuous  curve  $(u(r), \mu _1)(\xi _1)$ parameterized by $\xi _1 \in R$. If, in addition, $\lim _{|u| \ra \infty} \frac{g(uz)}{u}=0$ uniformly in $z \in R$, then $\frac{u(x)}{\xi _1} \ra \p _1(r)$ in $C^{2+\al }(B)$ as $\xi _1 \ra \pm \infty$. Moreover, all solutions of (\ref{30}) are radially symmetric, $u=u(r)$ with $r=|x|$, so that they satisfy
\beqa
\lbl{rad1}
&\s u''(r)+\frac{n-1}{r}u'(r)+\la _1 u+ g(u)=\mu _1 \p _1(r)+e(r) \,, \s \mbox{for $0<r<1$}  \\ \nonumber
& u'(0)=u(1)=0 \,,
\eeqa
\end{thm}

Except for the symmetry assertion this result was proved in \cite{K1}, where more general domains and non-radial  $e=e(x)$ were considered. Here radial symmetry follows from the uniqueness of the solution curve. Indeed, if a non-symmetric solution existed, any of its rotations would produce a different solution of (\ref{30})  with the same first harmonic $\xi _1$, hence lying on a different solution curve, but there is only one solution curve.
\medskip

We now discuss oscillations of the curve of radial solutions.
\medskip

\begin{thm}\lbl{thm:20}
In addition to the conditions of the Theorem \ref{thm:1}, assume that $g(u)$ is a periodic function of mean zero. Then $\mu _1(\xi _1) \ra 0$ as $|\xi _1| \ra \pm \infty$. The oscillation properties of the solution curve $\mu _1(\xi _1)$ of (\ref{thm:1}) depend on the dimension $n$ as follows.
\smallskip

\noindent
(i) For $1 \leq n \leq 3$, $\mu _1(\xi _1)$ changes sign infinitely many times on $(0,\infty)$.
 \smallskip

\noindent
(ii) In case $n=4$, $\mu _1(\xi _1)$ changes sign infinitely many times, provided that the same is true for the function $- \frac{1}{\p _1'(1)} g_2(0)+ \frac{2}{{\p _1 ''}^2(0)}  g_2 \left(\xi  \right)$ which includes the case $g_2(0)=0$), and there are only finitely many sign changes otherwise.
\smallskip

\noindent
(iii) For  $n=5$, in case $g_2(0) \ne 0$ there are only finitely many  sign changes of  $\mu _1(\xi _1)$, and if $g_2(0) = 0$, $\mu _1(\xi _1)$ changes sign infinitely many times. 
\smallskip

\noindent
((iv) For $n = 6$, in case $g_2(0) \ne 0$ there are only finitely many  sign changes of  $\mu _1(\xi)$, and if $g_2(0) = 0$,  $\mu _1(\xi)$ changes sign infinitely many times, provided that the same is true for the function $ f_3(1)g_3(0)-  f_3(0) g_3\left(\xi  \right)$ (which includes the case $g_3(0)=0$), and there are only finitely many sign changes otherwise.
\smallskip

\noindent
(v) For $n \geq 7$, assume that  $g_2(0) = 0$, but  $g_3(0) \ne 0$. Then the number of sign changes of   $\mu _1(\xi)$ is at most finite.
\end{thm}

\pf
We begin by sketching the proof. Since solutions of (\ref{30}) are radially symmetric, and $v(r)=\frac{u(r)}{\xi _1} \ra \p _1(r)$ in $C^{2+\al }(B)$, it follows that $u(r)$ is unimodular for large $\xi _1$, with a global maximum at $r=0$. It also follows that derivatives of $v(r)$, up to the order four,  tend to the corresponding derivatives of $\p _1(r)$ as $\xi _1 \ra \infty$ (expressing the derivatives of order greater than two from the equation for $v(r)$). Then proceed as in the Theorem \ref{thm:2}.
\medskip

Let $\eta =u(0)$, the maximum value of $u(r)$  for large $\xi _1$, and set $u(r)=\eta v(r)$, so that $v(0)=1$,  $\frac{\eta}{\xi _1} \ra 1$ and $v(r) =\p _1(r)+o(1)$ as $\xi _1 \ra \infty$. From (\ref{30}) obtain
\beq
\lbl{30a}
\s\s\s\; v''+\frac{n-1}{r}v'+\la _1 v +\frac{1}{\eta} g(\eta v)= \frac{\mu _1}{\eta} \p _1+\frac{1}{\eta}e  \s \mbox{for $x \in B$}, \s u=0 \s \mbox{on $\partial B$} \,.
\eeq
We now study the oscillations of $\mu _1=\mu _1(\xi _1)$ as $\xi _1 \ra \pm \infty$. Multiplying  the PDE version of the equation (\ref{30a}) by $\p _1(r)$ and integrating over the ball $B$ gives
\beq
\lbl{31}
\mu _1(\eta)=\frac{\int_0^1 g(\eta v(r)) \p _1(r) r^{n-1} \, dr}{\int_0^1  \p^2 _1(r) r^{n-1} \, dr} \,.
\eeq

The number of oscillations of the integral $\int_0^1 g(\eta v(r)) \p _1(r) r^{n-1} \, dr$ will depend on the dimension $n$. We  proceed as in the Theorem \ref{thm:2}, using up to three integrations by parts (the formulas (\ref{2}), (\ref{3}) and (\ref{4})), depending on the dimension $n$. Here $f_1(r)=\frac{r^{n-1} \p _1(r)}{v'(r)}$, $f_2(r)=\frac{f_1'(r)}{v'(r)}$ and $f_3(r)=\frac{f_2'(r)}{v'(r)}$. In the Theorem \ref{thm:2} we had instead: $g_1(r)=\frac{r^{n-1} \p _1(r)}{\p _1'(r)}$, $g_2(r)=\frac{g_1'(r)}{\p _1'(r)}$ and $g_3(r)=\frac{g_2'(r)}{\p _1'(r)}$ (we changed the notation of these functions to avoid confusion between the old and the new $f_i(r)$). As in \cite{G} we show that for $n \geq 2$ the function $f_1(r)$ is of class $C^{\infty}$,  $f_2(r) \in C^{\infty}$ for $n \geq 4$ , and $f_3(r) \in C^{\infty}$  for $n \geq 6$. Indeed, the proof in  \cite{G} was using only that $v(r) \ra \p _1(r)$ as $\eta \ra \infty$, which is true here too.
\medskip

Denote $F(v(r),r)=\frac{\mu _1}{\eta} \p _1(r)+\frac{1}{\eta}e(r)-\frac{1}{\eta} g(\eta v(r))$. As $\eta \ra \infty$,  $F(v(r),r) \ra 0$, and then $v(r) \ra \p _1(r)$ in $C^{2+\al }(B)$, as was pointed out previously. Expressing higher derivatives of $v(r)$ from (\ref{30a}), we see that the third and the fourth derivatives  of $v(r)$ at $r=0$ and at $r=1$ tend to the corresponding values of $\p _1(r)$. It follows that  at $r=0$ and $r=1$ the functions $f_1(r), f_2(r), f_3(r)$ tend to the corresponding values of  $g_1(r), g_2(r), g_3(r)$. Then all of the conclusions are the same as in the Theorem \ref{thm:2}.
\epf

\noindent
{\bf Remark } We now elaborate on the last step of the proof above. For the pivotal case $n=5$ we show directly that $g'_2(0)$ tends to a non-zero quantity as $\eta \ra \infty$, so that the argument proceeds as in Theorem \ref{thm:2}. Write (\ref{30a}) as
\beq
\lbl{33}
v''(r)+\frac{n-1}{r}v'(r) +f(r,v(r))=0 \,,
\eeq
where $f(r,v(r))=\la _1 v(r)+\frac{1}{\eta} g(\eta v(r))- \frac{\mu _1}{\eta} \p _1(r)-\frac{1}{\eta}e(r)$. Expressing $v''(r)$ from (\ref{33}), calculate
\[
f_2(r)=\frac{r^{n-1}f(r,v(r)) \p _1(r)+2(n-1)r^{n-2}\p _1(r) v'(r)+r^{n-1} \p '_1(r)v'(r)}{{v'}^3(r)} \,,
\]
\[
f'_2(r)=\frac{3r^{n-1}f^2(r,v(r)) \p _1(r)+8(n-1)r^{n-2}f(r,v(r)) \p _1(r)v'(r)+3 r^{n-1}f \p '_1v' }{{v'}^4(r)}
\]
\[
+\frac{r^{n-3} \left(8-14n+6n^2-r^2 \la _1\right)\p _1 v'+r^{n-1}  \p _1(r)\frac{d}{dr}f +4(n-1) r^{n-2} \p' _1(r)v'(r)  }{{v'}^3(r)} \,.
\]
Then for $n=5$
\[
f'_2(0)=\frac{1}{{v''}^4(0)} \left[  3f^2(0,1)  +32f(0,1)v''(0)+\al {v''}^2(0)  \right] \,,
\]
where $\al =6n^2-14n+8 \, {\large |}_{_{n=5}}=88$. In case $n=5$, $v''(0)=-\frac{1}{5} f(0,1)$, giving 
\[
f'_2(0)=\frac{75}{f^2(0,1)} \,,
\]
and $f(0,1)=\la _1 +\frac{1}{\eta} g(\eta )- \frac{\mu _1}{\eta} -\frac{1}{\eta}e(0) \ra \la _1$ as $\eta \ra \infty$. One can proceed similarly for other dimensions $n$.

\section{A more general result in one dimension}
\setcounter{equation}{0}
\setcounter{thm}{0}
\setcounter{lma}{0}

We  begin with a more general result for the oscillating integral
\beq
\lbl{1d1}
I(\xi )=\int _0^{\pi} g(\xi \sin x) \sin x \, dx \,,
\eeq
where $g(u)$ is not assumed to be periodic. 
\begin{thm}\lbl{thm:5}
Assume that $g(u) \in C(R)$   has finitely many roots on any bounded interval, and it changes sign at each root. Define $H(u)=\int _0^u g(t) t \, dt$. Assume that there exist two sequences $\{ \xi _n\} \ra \infty$ and $\{ \eta _n\} \ra \infty$ such that $H(\xi _n)>0$ and $H(\eta _n)<0$. Then $I(\xi )$ changes sign infinitely many times as $\xi \ra \infty$.
\end{thm}

\pf
Clearly, $g(u)$ has infinitely many roots. In case $g(0)=0$, let  $u_1 \geq 0$ be the supremum of $\beta$'s so that $g(u) \equiv 0$  on $[0,\beta)$, and in case $g(0) \ne 0$, set $u_1=0$. Clearly,  there exists $u_2>u_1$ so that $g(u)$ is either positive or negative on $(u_1,u_2)$. Without loss of generality we can make the following three assumptions.
\medskip

\noindent
(a) $g(u)$ is  negative on $(u_1,u_2)$. Otherwise consider $-I(\xi)$ involving $-g(u)$.
\medskip

\noindent
(b) $g(\xi _n) \geq 0$, and if $g(\xi _n) = 0$ then $g(\xi)$ is positive to the left of $\xi _n$. Otherwise let $\bar \xi _n$ be the largest root of $g(u)$ to the left of $\xi _n$. Then $H(\bar \xi _n)>H( \xi _n)>0$. Replace $\xi _n$ by $\bar \xi _n$.
\medskip

\noindent
(c)  $g(\eta _n) \leq 0$, and if $g(\eta _n) = 0$, then $g(\xi)$ is negative  to the left of $\xi _n$. The justification is similar to that of part (b).
\medskip

Setting $y=\sin x$ on the intervals $(0,\frac{\pi}{2})$ and $(\frac{\pi}{2},\pi)$, express
\beq
\lbl{1d2}
I(\xi )=2\int _0^1 g(\xi y) y \frac{1}{\sqrt{1-y^2}} \, dy \,.
\eeq

We claim that $I(\xi _n)>0$. Observe first that
\[
\int _0^1 g(\xi _ny) y  \, dy=\frac{1}{\xi _n^2} H(\xi _n)>0 \,.
\]
The graph of the function $g(\xi _n y) y$ on $(0,1)$ consists of pairs of humps. Each pair consists of a negative hump followed by positive hump. The function $\frac{1}{\sqrt{1-y^2}}>1$ is increasing, and so it  favors the positive humps. It follows that
\[
I(\xi _n)>2\int _0^1 g(\xi _ny) y  \, dy>0 \,.
\]

We claim that $I(\eta _n)<0$. Again, we begin by observing that
\[
\int _0^1 g(\eta _ny) y  \, dy=\frac{1}{\eta _n^2} H(\eta _n)<0 \,.
\]
The graph of the function $g(\eta _n y) y$ begins with a negative hump, and then it has pairs of positive humps followed by negative humps. The function $\frac{1}{\sqrt{1-y^2}}>1$ is increasing, and so it  favors the negative humps in each pair, while over the first negative hump the function $\frac{1}{\sqrt{1-y^2}}>1$ makes the integral smaller. It follows that
\[
I(\eta _n)<2\int _0^1 g(\eta _ny) y  \, dy<0 \,,
\]
completing the proof.
\epf

We now consider the problem (\ref{30}) in one dimension, which is convenient to consider on the interval $(0,\pi)$
\beq
\lbl{1d3}
\s\s u''+u+g(u)=\mu _1 \sin x+e(x) \,, \s \mbox{for $x \in (0,\pi)$}, \s u(0)=u(\pi)=0 \,.
\eeq
The problem is at resonance. Here $\la _1=1$, $\p _1(x)=\sin x$, $\la _2=4$. We assume that the function $e(x) \in C^{\alpha}$, $\alpha >0$, is even with respect to $\frac{\pi}{2}$ and satisfying  $\int _0^{\pi} e(x) \sin x \, dx=0$.
As above, decompose $u(x)=\xi _1 \sin x +U(x)$, with $\int _0^{\pi} U(x) \sin x \, dx=0$.
\medskip

Applied to (\ref{1d3}) the Theorem \ref{thm:1} implies the following result.
\begin{thm}\lbl{thm:6}
Assume that $g(u) \in C^2(R)$ satisfies the conditions (\ref{2a1}) and  (\ref{2a2}), with $\la _1=1$ and $\la _2=4$, and also  that $\lim _{|u| \ra \infty} \frac{g(uz)}{u}=0$ uniformly in $z \in R$.
Then the solution set of (\ref{1d3}) consists of a single continuous  curve  $(u(x), \mu _1)(\xi _1)$ parameterized by $\xi _1 \in R$. Moreover, all solutions of (\ref{1d3}) are even functions with respect to $\frac{\pi}{2}$, and $\frac{u(x)}{\xi _1} \ra \sin x$ in $C^2(0,\pi)$ as $\xi _1 \ra \pm \infty$. 
\end{thm}

The solution curve of (\ref{1d3}) performs infinitely many oscillations around the origin under the following conditions.

\begin{thm}\lbl{thm:7}
Assume that $g(u) \in C(R)$   has finitely many roots on any bounded interval, and it changes sign at each root.
Denoting $H(u)=\int _0^u g(t) t \, dt$, assume in addition to the conditions of the Theorem \ref{thm:6} that there exist two sequences $\{ \xi _n\} \ra \infty$ and $\{ \eta _n\} \ra \infty$ such that $H(\xi _n)>\ep$ and $H(\eta _n)<-\ep$, for some $\ep >0$. Then the function $\mu _1=\mu _1 (\xi _1)$ changes sign infinitely many times, as $\xi _1 \ra \infty$. In particular, at $\mu _1 =0$, the problem (\ref{1d3}) has infinitely many solutions.
\end{thm}

\pf
In view of the Theorem \ref{thm:6} only the last statement needs to be proved. Since  $\frac{u(x)}{\xi _1} \ra \sin x$ in $C^2(0,\pi)$ as $\xi _1 \ra \pm \infty$, and   $u(x)$ is  even with respect to $\frac{\pi}{2}$,  it follows that $u(x)$ is unimodular with a point of global maximum at  $\frac{\pi}{2}$, for large $\xi _1$. Let $\eta$ denote the maximum value of $u(x)$, and set $u(x)=\eta v(x)$ in (\ref{1d3}) to obtain
\beq
\lbl{1d5-}
\s\s v''+v+\frac{1}{\eta}g( \eta v)=\frac{\mu _1}{\eta} \sin x+\frac{1}{\eta} e(x) \,, \s \mbox{$x \in (0,\pi)$}, \s v(0)=v(\pi)=0 \,.
\eeq
Clearly $v(\frac{\pi}{2})=1$, $\frac{\eta}{\xi _1}  \ra 1$ and $v(x)=\sin x+o(1)$, as $\xi _1 \ra \infty$. Multiplication of (\ref{1d5-}) by $\sin x$ and integration over $(0,\pi)$ gives
\beqa \nonumber
& \mu _1 \frac{\pi}{2}=\int _0^{\pi} g( \eta v(x)) \sin x \, dx =2\left(1+o(1) \right) \int _0^{\frac{\pi}{2}} g( \eta v(x)) v(x) \, dx\\ \nonumber
& =2\left(1+o(1) \right) \int _0^1 g( \eta y) y \frac{dx}{dy}\, dy \,,
\eeqa
setting $y=v(x)$. The function $\frac{dx}{dy}$ tends to a positive increasing function that is greater than $1$ for $y \in (0,1)$ (namely, to $\frac{1}{\sqrt{1-y^2}}$). As in the Theorem \ref{thm:5}, $\mu _1$ changes sign infinitely many times as $\eta \ra \infty$.
\epf

\medskip

\noindent
{\bf Example 1 }  We computed the solution curve $\mu _1= \mu _1(\xi _1)$, $\xi _1>0$,  for the following example, with the linear part at resonance,
\beq
\lbl{1d5}
\s\s\s\s u''+u+\frac{\sin u}{\sqrt{u+4}}=\mu _1 \sin x+\sin 3x \,, \s \mbox{for $x \in (0,\pi)$}, \s u(0)=u(\pi)=0 \,.
\eeq
(Recall that the full solution curve of (\ref{1d5}) has the form $\left(u(x),\mu _1 \right)(\xi _1)$, where $u(x)=\xi _1 \p _1 (x)+U(x)$.)
Here $\la _1=1$, $\p _1(x)=\sin x$. It is easy to check that the Theorem \ref{thm:7} applies, so that there are infinitely many solutions at $\mu _1=0$. The oscillating solution curve $\mu _1= \mu _1(\xi _1)$ is presented in Figure \ref{fig:1}, see \cite{KS1} for the listing of the {\em Mathematica} program used, together with detailed explanations.

\begin{figure}
\begin{center}
\scalebox{0.99}{\includegraphics{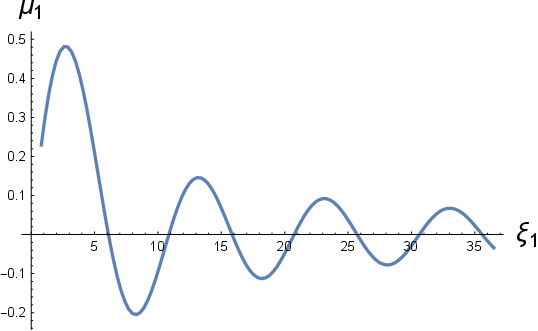}}
\end{center}
\caption{ The  solution curve $\mu _1= \mu _1(\xi _1)$ of the problem (\ref{1d5}), oscillating around the $\xi _1$-axis.}
\lbl{fig:1}
\end{figure}

\medskip

\noindent
{\bf Example 2 }  Theorem \ref{thm:7} does not apply to the problem
\beq
\lbl{1d6}
\s\s\s\; u''+u+\frac{\sin u}{\sqrt{u^4+4}}=\mu _1 \sin x+\sin 3x \,, \s \mbox{for $x \in (0,\pi)$}, \s u(0)=u(\pi)=0 \,.
\eeq
(Here the integral $\int_0^{\infty} g(t)t \, dt$ converges.) Our calculations, presented in Figure \ref{fig:3}, suggest that when $\xi _1>0$ there are no solutions for  $\mu _1=0$, while there are arbitrary many solutions for $\mu _1>0$ sufficiently small. So that if $g(u)$ in (\ref{1d3}) tends to zero sufficiently fast as $u \ra \infty$, the solution curve may tend to zero without oscillating around the origin.

\begin{figure}
\begin{center}
\scalebox{0.99}{\includegraphics{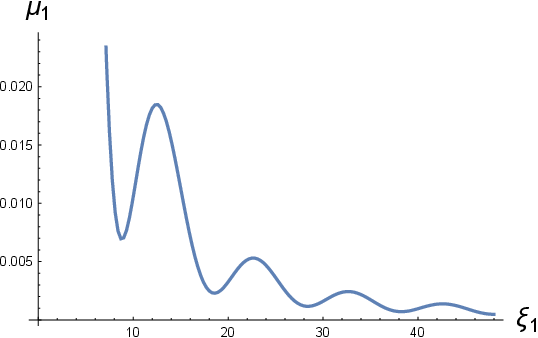}}
\end{center}
\caption{ The  solution curve $\mu _1= \mu _1(\xi _1)$ of the problem (\ref{1d6}), oscillating near $\mu _1=0$}
\lbl{fig:3}
\end{figure}


\begin{thebibliography}{99}

\bibitem{C}
 D. Costa,  H. Jeggle,  R. Schaaf  and  K. Schmitt,     Oscillatory perturbations of linear     
      problems at resonance, {\em Results in Mathematics} {\bf 14},  275-287  (1988).

\bibitem{C0}
O. Costin, N. Falkner and J.D.  McNeal, 
Some generalizations of the Riemann-Lebesgue lemma,
{\em Amer. Math. Monthly} {\bf  123}, no. 4, 387-391  (2016).


\bibitem{D1} 
E.N. Dancer, On the use of asymptotics in nonlinear boundary value problems, {\em  Ann. Mat. Pura Appl.} {\bf 131}, (4), 167-185   (1982).


\bibitem{G}
A. Galstian, P.   Korman   and   Y. Li,     On the oscillations of the solution curve for a class of semilinear equations, {\em   J. Math. Anal. Appl.} {\bf  321},  no. 2,  576-588 (2006).

\bibitem{ks}
K. Kato and T. Shibata,  Simple proof of stationary phase method and application to oscillatory bifurcation problems, {\em Nonlinear Anal.} {\bf  190} (2020), 111594, 13 pp.

\bibitem{K}
P. Korman, Global Solution Curves for Semilinear Elliptic Equations, World Scientific, Hackensack, NJ (2012).

\bibitem{K1}
P. Korman, Global solution curves in harmonic parameters, and multiplicity of solutions, {\em J. Differential Equations} {\bf 296}, 186-212  (2021).

\bibitem{KS}
P. Korman and D.S. Schmidt, Infinitely many solutions and asymptotics for resonant oscillatory problems, Special issue in honor of Alan C. Lazer, Electron. J. Diff. Equ., Special Issue 01, 301-313 (2021).

\bibitem{KS1}
P. Korman and D.S. Schmidt, Calculating global solution curves for boundary value problems, {\em Wolfram Notebook Archive}, online: notebookarchive.org/calculating-global-solution-curves-for-boundary-value-problems--2022-08-eb98nqk/.

\bibitem{L}
E.M. Landesman and A.C. Lazer,   Nonlinear perturbations of linear elliptic        
 boundary value problems at resonance, {\em J. Math. Mech.} {\bf 19}, 609-623 (1970).

\bibitem{leach}
A.C. Lazer  and  D.E. Leach,  Bounded perturbations of forced harmonic oscillators at resonance, {\em Ann. Mat. Pura Appl.}  {\bf 82} (4),  49-68 (1969). 

\bibitem{SS}
R. Schaaf  and K. Schmitt, A class of nonlinear Sturm-Liouville   
      problems with infinitely many solutions, {\em Trans. Amer. Math. Soc.} {\bf 306}, no. 2,  853-859 (1988).


\bibitem{SS1}
R. Schaaf  and K. Schmitt, Asymptotic behavior of positive solution branches of elliptic problems with linear part at resonance, {\em Z. Angew. Math. Phys.} {\bf  43}, no. 4, 645-676  (1992).

\bibitem{sh}
T. Shibata,  Asymptotic expansion of oscillatory bifurcation curves of ODEs with nonlinear diffusion, {\em Differential Integral Equations} {\bf 33}, no. 5-6, 257-272  (2020).

\end{thebibliography}
\end{document}